\documentclass{amsart}

\newtheorem{theorem}{Theorem}
\newcommand{\bt}{\begin{theorem}}
\newcommand{\et}{\end{theorem}}

\newtheorem{lemma}{Lemma}
\newcommand{\bl}{\begin{lemma}}
\newcommand{\el}{\end{lemma}}

\newtheorem{corollary}{Corollary}
\newcommand{\bc}{\begin{corollary}}
\newcommand{\ec}{\end{corollary}}

\newtheorem{problem}{Problem}
\newcommand{\bprob}{\begin{problem}}
\newcommand{\eprob}{\end{problem}}

\DeclareMathOperator{\card}{\text{card}}

\DeclareMathOperator{\qqand}{\qquad\text{and}\qquad}

\newcommand{\E}{\ensuremath{\mathbf E }}

\newcommand{\Z}{\ensuremath{\mathbf Z}}

\newcommand{\mcg}{\ensuremath{ \mathcal G}}
\newcommand{\mcu}{\ensuremath{ \mathcal U}}

\newcommand{\beq}{\begin{equation}}
\newcommand{\eeq}{\end{equation}}
\newcommand{\benum}{\begin{enumerate}}
\newcommand{\eenum}{\end{enumerate}}

\DeclareMathOperator{\Prob}{\mathbf{Pr}}

\usepackage{amsmath,amssymb,amsthm}

\title[Integers with missing digits]{Convergent series of integers with missing digits}
\author{Melvyn B. Nathanson}
\address{Lehman College (CUNY), Bronx, NY 10468}\email{melvyn.nathanson@lehman.cuny.edu}



\date{\today}

\begin{document}

\begin{abstract}
A classical theorem of Kempner states that the sum of the reciprocals of positive integers 
with missing decimal digits converges.  
This result is extended to much larger families of ``missing digits'' sets of positive integers  
with convergent  harmonic series.  
\end{abstract}

\maketitle 

\section{Kempner's theorem for $g$-adic repesentations}

A \emph{harmonic series} is a series of the form $\sum_{a\in A} 1/a$, 
where $A$ is a set of positive integers.  
For example, let $A_{10,9}$ be the set of positive integers 
in which the digit $9$ does not occur in the usual decimal representation.  
Kempner~\cite{kemp14}  proved in 1914 that the harmonic series 
$\sum_{a\in A_{10,9}} 1/a$ converges.  
He called this  ``a curious convergent series.''  
More generally, for every integer $g \geq 2$,  every positive integer $n$ 
has a unique  \textit{$g$-adic representation} of the form $n = \sum_{i=0}^{m-1} c_i g^i$, 
with digits $c_i \in \{0,1,2,\ldots, g-1\}$ for $i = 0,1,\ldots, m-1$ and $c_{m-1} \neq 0$. 
If $A_{g,c}$ is the set of integers whose $g$-adic representation contains no digit $c$, then 
the infinite series $\sum_{a\in A_{g,c}} 1/a$ converges.  
This is Theorem 144 in Hardy and Wright~\cite{hard-wrig08}.

Kempner's theorem has been studied and extended by 
Baillie~\cite{bail79}, Farhi~\cite{fahr08}, Gordon~\cite{gord19}, Irwin~\cite{irwi16}, 
Nathanson~\cite{nath21y}, 
Schmelzer and Baillie~\cite{schm-bail08}, and Wadhwa~\cite{wadh75,wadh78}.  

Lubeck and Ponomarenko~\cite{lube-pono18} used a probability argument to obtain 
a remarkably strong form of Kempner's theorem. 

\bt[Lubeck and Ponomarenko]           \label{MissingDigits:theorem:LP} 
Let $0 < \lambda < 1$. 
Let $A^{(\lambda)}_m$ be the set of integers $n \in  [10^{m-1}, 10^m-1]$ 
such that at most $\lambda m$ 
digits of $n$ are equal to 9, and let $A_{10,9}^{(\lambda)} = \bigcup_{m=1}^{\infty} A^{(\lambda)}_m$.  
If $\lambda < 1/10$, then the harmonic series $\sum_{a\in A_{10,9}^{(\lambda)}} 1/a$ converges. 
If $\lambda \geq 1/10$, then the harmonic series $\sum_{a\in A_{10,9}^{(\lambda)}} 1/a$ diverges. 
\et

This result extends to $g$-adic representations for all $g \geq 2$.  

The $g$-adic representation is a special case of a more general method to represent the positive integers.  
A \emph{\mcg-adic sequence} is a strictly increasing sequence of positive integers $\mcg = (g_i)_{i=0}^{\infty}$ 
such that  $g_ 0 = 1$ and $g_{i}$ divides $g_{i+1}$ for all $ i \geq 0$.   
The integer
\[
d_{i} = \frac{g_{i+1}}{g_{i}}
\]
satisfies $d_i \geq 2$ and 
\beq               \label{MissingDigits:gk-product}
g_{k+1} = g_{k}d_{k} = d_0 d_1 d_2 \cdots d_{k-1}d_k
\eeq
for all $k \geq 0$.
Every positive integer $n$ has a unique representation in the form 
\beq               \label{MissingDigits:n}
n = \sum_{i=0}^{m-1} c_i g_i
\eeq
where $c_i \in \{0,1,\ldots, d_i-1\}$ for all $i  = 0,1,\ldots, m-1$ and $c_{m-1} \neq 0$. 
We call~\eqref{MissingDigits:n} the \textit{$\mcg$-adic representation} of $n$. 
This is equivalent to deBruijn's additive system (Nathanson~\cite{nath2014-150, nath2017-172}).

Harmonic series constructed from sets of positive integers 
with missing $\mcg$-adic digits do not necessarily  converge. 
The following result constructs sets of integers 
with missing $\mcg$-adic digits whose harmonic series do converge, 
and also sets of integers 
with missing $\mcg$-adic digits whose harmonic series do not diverge. 

The \mcg-adic sequence $\mcg = (g_i)_{i=0}^{\infty}$ has 
\emph{quotients bounded  by $d$} if  
\beq                \label{Digits:limsup}
d_i = \frac{g_{i+1}}{g_i}  \leq d\qquad \text{for all $i =0,1,2,\ldots$.} 
\eeq
Nathanson~\cite{nath21x} obtained the following criterion for the convergence and divergence 
of harmonic series of positive integers with missing \mcg-adic digits.

\begin{theorem}[Nathanson]                   \label{MissingDigits:theorem:converges} 
Let $\mcg = (g_i)_{i=0}^{\infty}$ be a \mcg-adic sequence with  quotients bounded by $d$.  
Let $I$ be a set of nonnegative integers, and 
let $I(k)$ be the number of $i \in I$ with $i \leq k$.  
For all $i \in I$, let $U_i$ be a nonempty proper subset of  $[0,d_i-1]$.  

Let $n = \sum_{i=0}^k c_i g_i$ be the $\mcg$-adic representation of the positive integer $n$, 
and let $A$ be the set of positive integers $n$ such that 
$c_i \in [0,d_i -1]  \setminus U_i$ for all $i \in I$.  
If 
\beq                \label{Digits:I(x)-1}
I(k) \geq \frac{(1+\delta)\log k}{\log (d/(d-1))}
\eeq 
for some $\delta > 0$ and all $k \geq k_0(\delta) $, 
then  the harmonic series $\sum_{a\in A} 1/a$ converges.
If 
\beq                \label{Digits:I(x)-2}
I(k) \leq \frac{(1-\delta) \log k}{\log d}
\eeq 
for some $\delta> 0$ and all $k \geq k_1(\delta) $, 
then the harmonic series $\sum_{a\in A} 1/a$ diverges.
\end{theorem}

In this paper we show that the Lubeck-Ponomarenko theorem extends to \mcg-adic representations.    
Moreover, it can also be significantly strengthened, even for the usual $g$-adic representations.

\section{Harmonic series for \mcg-adic repesentations} 

We use a simple form of Chernoff's inequality in probability theory.

\bt[Chernoff]
Let $p_0, p_1,\ldots, p_{m-1} \in [0,1]$  and let $X_0, X_1,\ldots, X_{m-1}$ 
be independent Bernoulli random variables with  
\[
\Prob(X_i =1) = p_i \qqand \Prob(X_i=0) = 1-p_i  
\]
for $i = 0, 1, \ldots, m-1$.  Let 
\[
X = X_0 + X_1 + \cdots + X_{m-1} 
\] 
and 
\[
\mu = \E(X) = p_0 + p_1+\cdots + p_{m-1}. 
\]
For all $0 < \delta < 1$, 
\[
\Prob(X \leq (1-\delta) \mu) < e^{-\mu \delta^2/2}.   
\]
\et

The short proof is in the Appendix.

We denote the cardinality of the set $S$ by $\card(S)$ or by $|S|$.  

\bt                    \label{MissingDigits:theorem:A-lambda} 
Let $\mcg = (g_i)_{i=0}^{\infty}$ be a  \mcg-adic sequence 
with quotients bounded  by $d$, and let 
\beq               \label{MissingDigits:lambda}
0 \leq \lambda < \frac{1}{d}. 
\eeq
Let $U_i$ be a proper nonempty subset of $[0,d_i-1]$ for all $i = 0,1,2,\ldots$, 
and let $\mcu = (U_i)_{i=0}^{\infty}$.  
Let  $A_{\mcg,\mcu}^{(\lambda)}$ be the set of  positive integers $n$ such that, 
if $n \in  [g_{m-1},g_m-1]$ has the \mcg-adic representation 
\[
n = \sum_{i=0}^{m-1} c_i g_i 
\]
then
\[
\card\{i\in [0,m-1]: c_i \in U_i \} \leq \lambda m.
\]
There exists $\sigma < 1$ such that 
\[
\sum_{a \in A_{\mcg,\mcu}^{(\lambda)} } \frac{1}{a^{\sigma}} < \infty.
\]
\et

\begin{proof}
Let 
\[
I_m = [g_{m-1},g_m-1] \cap \Z \qqand J_m = [0,g_m-1] \cap \Z.
\]
We have $|I_m| = g_m-g_{m-1}$ and $|J_m| = g_m$.  
Every integer $n \in J_m$ has a unique \mcg-adic representation 
\[
n = \sum_{i=0}^{m-1} c_i g_i 
\]
with $c_i \in [0,d_i-1]$ for all $i = 0,1,\ldots, m-1$.   
Note that $I_m \subseteq J_m$ and $n \in J_m \setminus I_m$ if and only if $c_{m-1} = 0$.

Define the uniform probability measure on the set $J_m$: 
\[
\Pr(n) = \frac{1}{|J_m|} = \frac{1}{g_m} \qquad \text{for all $n \in J_m$.}   
\]
Define the random variable $X:J_m \rightarrow [0,m]$ as follows.  
For $n = \sum_{i=0}^{m-1} c_i g_i  \in J_m$, let 
\beq               \label{MissingDigits:X}
X(n) = \card\{i\in [0,m-1]: c_i \in U_i \}. 
\eeq
We have 
\[
\Pr \left( X \leq \lambda m \right) = \frac{\left| A_{\mcg,\mcu}^{(\lambda)} \cap J_m \right| }{|J_m|} 
\]
and so 
\beq               \label{MissingDigits:AmBm} 
\left| A_{\mcg,\mcu}^{(\lambda)} \cap I_m \right| 
 \leq \left| A_{\mcg,\mcu}^{(\lambda)} \cap J_m \right|  
= \Pr \left( X \leq \lambda m \right) g_m.
\eeq

For all $i = 0,1,2,\ldots,m-1$, define the Bernoulli random variable $X_i: J_m \rightarrow \{0,1\}$ by 
\beq               \label{MissingDigits:Xi}
X_i (n) = \begin{cases}
1 & \text{if $c_i \in U_i$} \\ 
0 & \text{if $c_i \notin U_i$.}
\end{cases}
\eeq
We have 
\[
\Pr( X_i = 1) = \frac{|U_i|}{d_i}  
\]
and 
\[
\Pr( X_i = 0) = 1 -  \frac{ |U_i|}{d_i} 
\]
and so $X_i$ has the expected value 
\[
\E(X_i) = \frac{|U_i|}{d_i}. 
\]
Choose $\varepsilon_i \in \{0,1\}$ for all $i \in \{0,1,\ldots,m-1\}$.  
Recall that 
\[
|J_m| = g_m = \prod_{i=0}^{m-1} d_i.
\]
We have 
\begin{align*} 
\Pr ( X_i = \varepsilon_i  & \text{ for all } i \in \{0,1,\ldots,m-1\} ) \\
& = \frac{1}{|J_m|} \prod_{\substack{i =0 \\ \varepsilon_i = 1}}^{m-1} |U_i|
\ \prod_{\substack{i =0 \\ \varepsilon_i = 0}}^{m-1} \left( d_i - |U_i|  \right)  \\
& = \prod_{\substack{i =0 \\ \varepsilon_i = 1}}^{m-1} \frac{|U_i|}{d_i} 
\ \prod_{\substack{i =0 \\ \varepsilon_i = 0}}^{m-1} \left(  1 -  \frac{ |U_i|}{d_i}  \right)  \\
& = \prod_{{i =0 }}^{m-1} \Pr \left( X_i = \varepsilon_i  \right) 
\end{align*} 
and so $X_0, X_1, \ldots, X_{m-1}$ are independent random variables.

For $m = 1,2,3,\ldots$, define $d_m^*$ by 
\[
\frac{1}{d_m^*} = \frac{1}{m} \sum_{i=0}^{m-1} \frac{|U_i|}{d_i} \geq \frac{1}{m} \sum_{i=0}^{m-1} \frac{1}{d_i} 
\geq \frac{1}{m} \sum_{i=0}^{m-1} \frac{1}{d} =\frac{1}{d}.
\]
It follows that  $0 < d^*_m \leq d$.  
From~\eqref{MissingDigits:X} and~\eqref{MissingDigits:Xi} we obtain 
\[
X = \sum_{i=0}^{m-1} X_i 
\]
and so 
\[
\mu = \E(X)= \sum_{i=0}^{m-1} \E(X_i) = \sum_{i=0}^{m-1}\frac{|U_i|}{d_i} = \frac{m}{d^*_m} \geq \frac{m}{d}. 
\]
Because $X_0, X_1, \ldots, X_{m-1}$ are independent Bernoulli random variables, 
 Chernoff's inequality will give an upper bound for $\Pr(X \leq \lambda m)$.  
From~\eqref{MissingDigits:lambda}, we have 
\[
0 < \lambda d_m^* \leq \lambda d < 1
\]
and so 
\[
0  <  1 - \lambda d_m^* = \delta <  1.
\]
Thus, 
\[
\lambda m = \lambda d_m^* \mu = (1-\delta) \mu
\]
and  
\begin{align*}
\Pr(X \leq \lambda m) & = \Pr( X \leq (1-\delta) \mu ) \\ 
& < e^{-\mu\delta^2/2} = e^{-(1 - \lambda d_m^* )^2 m/2d_m^*} \\ 
& \leq e^{-(1 - \lambda d)^2 m/2d}.  
\end{align*} 
Inequality~\eqref{MissingDigits:AmBm} gives 
\beq                                 \label{MissingDigits:ineq-A} 
\left| A_{\mcg,\mcu}^{(\lambda)}\cap I_m \right| 
 \leq \Pr \left( X \leq \lambda m \right) g_m 
\leq e^{-(1 - \lambda d)^2 m/2d} g_m.
\eeq
Choose  $\sigma \in (0,1)$ such that 
\beq                                 \label{MissingDigits:ineq-de} 
1 < d^{ 1-\sigma} < e^{(1 - \lambda d)^2 /2d}.
\eeq
Because \mcg\ has quotients bounded by $d$, we have 
\[
\frac{g_m}{d}  \leq \frac{g_m}{d_{m-1}} = g_{m-1} = \prod_{i=0}^{m-2} d_i \leq d^{m-1}. 
\]
Also, $a \geq g_{m-1}$ for all $a \in I_m$.  
Applying~\eqref{MissingDigits:ineq-A}, we obtain 
\begin{align*}
\sum_{ a\in A_{\mcg,\mcu}^{(\lambda)}}  \frac{1}{a^{\sigma} }
&=  \sum_{m=1}^{\infty} \sum_{a \in A_{\mcg,\mcu}^{(\lambda)}\cap I_m  }  \frac{1}{a^{\sigma} } \\
& \leq  \sum_{m=1}^{\infty} \frac{| A_{\mcg,\mcu}^{(\lambda)} \cap I_m   |}{g_{m-1}^{\sigma} }\\ 
& \leq  \sum_{m=1}^{\infty} \frac{ e^{-(1 - \lambda d )^2 m/2d} g_m }{g_{m-1}^{\sigma} }\\
& \leq d \sum_{m=1}^{\infty} \frac{ e^{-(1 - \lambda d )^2 m/2d} g_{m-1} }{g_{m-1}^{\sigma} }\\
& = d   \sum_{m=1}^{\infty} e^{-(1 - \lambda d )^2 m/2d} g_{m-1}^{1-\sigma} \\
& \leq d \sum_{m=1}^{\infty} e^{-(1 - \lambda d )^2 m/2d} d^{ (m-1)(1-\sigma)} \\
& =  d^{\sigma}   \sum_{m=1}^{\infty} \left( \frac{ d^{ 1-\sigma} }{ e^{(1 - \lambda d )^2 /2d} }\right)^m  \\
& < \infty
\end{align*}  
by inequality~\eqref{MissingDigits:ineq-de}.  This completes the proof.  
\end{proof}

\bc                             \label{MissingDigits:corollary:A-lambda} 
Let $A_{\mcg,\mcu}^{(\lambda)}$ be the set of positive integers constructed 
in Theorem~\ref{MissingDigits:theorem:A-lambda}.  The Dirichlet series 
\[
F_{A_{\mcg,\mcu}^{(\lambda)} }(s) = \sum_{a\in A_{\mcg,\mcu}^{(\lambda)} } \frac{1}{a^s} 
\]
has abscissa of convergence 
\[
\sigma_c \leq 1 - \frac{(1-\lambda d)^2 }{2d\log d}.
\]
\ec

Note that, in the $g$-adic case, Theorem~\ref{MissingDigits:theorem:LP}  
asserts only that $\sigma_c \leq 1$. 

\section{Open problems}

\bprob
Corollary~\ref{MissingDigits:corollary:A-lambda} 
gives an explicit upper bound for the abscissa of convergence 
of the Dirichlet series $F_{A_{\mcg,\mcu}^{(\lambda)} }(s)$.  Can this bound be improved?  
Can the abscissa of convergence be computed? 
\eprob

\bprob
Determine the convergence or divergence 
of the Dirichlet series $F_{A_{\mcg,\mcu}^{(\lambda)} }(s)$ for $\lambda \geq 1/d$.
\eprob

\bprob
Consider the convergence or divergence of Dirichlet series 
$F_{A_{\mcg,\mcu}^{(\lambda)} }(s)$  constructed from \mcg-adic sequences 
$\mcg = (g_i)_{i=0}^{\infty}$ with unbounded quotients $d_i = g_{i+1}/g_i$.  
\eprob

\bprob
Interpret Theorem~\ref{MissingDigits:theorem:A-lambda} in terms of the Hausdorff dimension 
of sets of real numbers with missing $g$-adic or \mcg-adic digits.
\eprob


\section*{Appendix: Proof of Chernoff's inequality}

Let $X$ be a nonnegative random variable on a countable probability space $\Omega$.  
For all $b > 0$, we have  Markov's inequality: 
\begin{align*}
\E(X) & = \sum_{\omega \in \Omega} X(\omega) \Prob(\omega) 
 \geq \sum_{\substack{\omega \in \Omega \\ X(\omega) \geq b}} X(\omega) \Prob(\omega) \\
& \geq b  \sum_{\substack{\omega \in \Omega \\ X(\omega) \geq b}} \Prob(\omega) 
  = b \Prob(X \geq b)   
\end{align*}
and so 
\beq                                  \label{MissingDigits:Markov} 
\Prob(X \geq b) \leq \frac{\E(X)}{b}.
\eeq 
For $s > 0$, the function $y = e^{-sx}$ is decreasing and positive, and so 
$X(\omega) \leq a$ if and only if $e^{-sX(\omega)} \geq e^{-sa}$.
Applying~\eqref{MissingDigits:Markov}  to the random variable $ e^{-sX} $ 
with $b = e^{-sa}$, we obtain 
\beq          \label{MissingDigits:Markov-2} 
\Prob(X \leq a) = \Prob\left( e^{-sX} \geq e^{-sa}  \right) 
\leq \frac{\E\left[e^{-sX} \right]}{e^{-sa}}.
\eeq

Let $X_i$ be a Bernoulli random variable with  
\[
\Prob(X_i =1) = p_i \qqand \Prob(X_i=0) = 1-p_i  
\]
Using the inequality $1 + x \leq e^{x}$, we have     
\[
\E\left[ e^{-sX_i} \right] = e^{-s}p_i+ e^0 (1-p_i) = 1 + \left( e^{-s} - 1 \right)p_i \leq e^{\left(e^{-s}-1\right)p_i}. 
\]

Let 
\[
X = X_0 + X_1 + \cdots + X_{m-1} 
\] 
and 
\[
\mu = \E(X) = p_0 + p_1+\cdots + p_{m-1}.
\]
Recall that if $Y_0, Y_1,\ldots, Y_{m-1}$ are independent random variables, then 
\[
\E\left[ \prod_{i=0}^{m-1} Y_i \right] =  \prod_{i=1}^{m-1} \E\left[ Y_i \right]. 
\]
The random variables $X_0, X_1,\ldots, X_{m-1}$ are independent, 
and so the random variables $e^{-sX_0}, e^{-sX_1},\ldots, e^{-sX_{m-1}} $ are independent.  Therefore, 
\begin{align*}
\E\left[ e^{-sX} \right] & = \E\left[ e^{-s\sum_{i=0}^{m-1} X_i} \right]
 =  \E\left[ \prod_{i=0}^{m-1} e^{-sX_i} \right] 
 =  \prod_{i =0}^{m-1} \E\left[  e^{-sX_i} \right]  \\ 
 & \leq  \prod_{i=0}^{m-1} e^{\left(e^{-s}-1\right)p_i} =  e^{\left(e^{-s}-1\right) \sum_{i=0}^{m-1}  p_i}  =  e^{\left(e^{-s}-1\right) \mu}. 
\end{align*}

Consider the function  
\[
h(\delta) = (1-\delta) \log(1-\delta) + \delta - \frac{\delta^2}{2}.
\]
We have  $h(0) = 0$ and, for $0 < \delta <1$, 
\[
h'(\delta) = -\log(1-\delta) - \delta > 0.  
\]
It follows that  $h(\delta) > 0$.  Equivalently, 
\beq                                  \label{MissingDigits:log-delta} 
-\delta - (1-\delta) \log(1-\delta)  < - \frac{\delta^2}{2}
\eeq
Let $s > 0$ and $0 < \delta < 1$.  
Inequality~\eqref{MissingDigits:Markov-2} with $a = (1-\delta)\mu$ gives 
\begin{align*}
\Prob(X \leq (1-\delta)\mu) 
&  \leq \frac{\E\left[e^{-sX} \right]}{e^{-s(1-\delta)\mu}} 
 \leq \frac{ e^{\left(e^{-s}-1\right) \mu} }{e^{-s(1-\delta)\mu}} \\
&  =  e^{\left( e^{-s}-1 + s(1-\delta) \right)\mu}.
\end{align*}
We minimize $e^{-s}-1 + s(1-\delta)$ by choosing  $s = - \log(1 - \delta)$.   
Applying  inequality~\eqref{MissingDigits:log-delta}, we obtain 
\begin{align*}
e^{-s}-1 + s(1-\delta) & = -\delta - (1 - \delta)\log(1 - \delta)  < - \frac{\delta^2}{2}  
\end{align*} 
and so 
\[
\Prob(X \leq (1-\delta)\mu) < e^{ -\mu\delta^2/2}. 
\]
This completes the proof.

\def\cprime{$'$} \def\cprime{$'$} \def\cprime{$'$}
\providecommand{\bysame}{\leavevmode\hbox to3em{\hrulefill}\thinspace}
\providecommand{\MR}{\relax\ifhmode\unskip\space\fi MR }
\providecommand{\MRhref}[2]{%
  \href{http://www.ams.org/mathscinet-getitem?mr=#1}{#2}
}
\providecommand{\href}[2]{#2}


\begin{thebibliography}{10}

\bibitem{bail79}
R. Baillie, \emph{Sums of reciprocals of integers missing a given digit},
  Amer. Math. Monthly \textbf{86} (1979), no.~5, 372--374.

\bibitem{fahr08}
B. Farhi, \emph{A curious result related to {K}empner's series}, Amer. Math.
  Monthly \textbf{115} (2008), no.~10, 933--938.

\bibitem{gord19}
R.~A. Gordon, \emph{Comments on ``{S}ubsums of the harmonic series''},
  Amer. Math. Monthly \textbf{126} (2019), no.~3, 275--279.

\bibitem{hard-wrig08}
G.~H. Hardy and E.~M. Wright, \emph{{An Introduction to the Theory of
  Numbers}}, 6th ed., Oxford University Press, Oxford, 2008.

\bibitem{irwi16}
F. Irwin, \emph{A curious convergent series}, Amer. Math. Monthly
  \textbf{23} (1916), no.~5, 149--152.

\bibitem{kemp14}
A.~J. Kempner, \emph{A curious convergent series}, Amer. Math. Monthly
  \textbf{21} (1914), no.~2, 48--50.

\bibitem{lube-pono18}
B. Lubeck and V. Ponomarenko, \emph{Subsums of the harmonic series},
  Amer. Math. Monthly \textbf{125} (2018), no.~4, 351--355.

\bibitem{nath2014-150}
M.~B. Nathanson, \emph{Additive systems and a theorem of de {B}ruijn}, Amer.
  Math. Monthly \textbf{121} (2014), no.~1, 5--17.

\bibitem{nath2017-172}
M.~B. Nathanson, \emph{Limits and decomposition of de {B}ruijn's additive
  systems}, Combinatorial and Additive Number Theory. {II}, Springer Proc.
  Math. Stat., vol. 220, Springer, 2017, pp.~255--267.



\bibitem{nath21y}
M.~B. Nathanson, \emph{Dirichlet series of integers with missing digits}, 
J. Number Theory, to appear.  


\bibitem{nath21x}
M.~B. Nathanson, \emph{Curious convergent series of integers with missing digits},
Integers, to appear.



\bibitem{schm-bail08}
T. Schmelzer and R. Baillie, \emph{Summing a curious, slowly convergent
  series}, Amer. Math. Monthly \textbf{115} (2008), no.~6, 525--540.

\bibitem{wadh75}
A.~D. Wadhwa, \emph{An interesting subseries of the harmonic series}, Amer.
  Math. Monthly \textbf{82} (1975), no.~9, 931--933.

\bibitem{wadh78}
A.~D. Wadhwa,\emph{Some convergent subseries of the harmonic series}, Amer. Math.
  Monthly \textbf{85} (1978), no.~8, 661--663.

\end{thebibliography}
\end{document}